# Discrete Distribution for a Wiener Process Range and its Properties

**Mohamed Abd Allah El-Hadidy** [1,2,*]

[1]Mathematics and Statistics Department, College of Science, Taibah University, Yanbu, Kingdom of Saudi Arabia

[2]Mathematics Department, Faculty of Science, Tanta University, Tanta, Egypt.

[*]**Corresponding author at:** Mathematics and Statistics Department, College of Science, Taibah University, Yanbu, Kingdom of Saudi Arabia.

**E-mail address:** melhadidi@science.tanta.edu.eg

**Abstract.** *We introduce the discrete distribution of a Wiener process range. Rather than finding some basic distributional properties including hazard rate function, moments, Stress-strength parameter and order statistics of this distribution, this work studies some basic properties of the truncated version of this distribution. The effectiveness of this distribution is established using a data set.*

## 1. Introduction

Discretizing continuous random variables have a remarkable importance in our life due to its great applicability. It is an important methodology for a lifetime modeling which may not necessarily always be measured on a continuous scale where the observed measurements are usually discrete in nature. Thus, the continuous lifetime may be often counted as a discrete random variable. One of the important real life examples is the difference between the highest price for the stock and it's the lowest price through days, months or years. Teamah et al. [1] presented the continuous distribution of this difference as the distribution of Wiener process range in bonded domain by using truncation method. Here, we need to obtain the distribution of a particular difference measured in days or months, etc.. Krishna and Singh [2] presented an interesting situation where the lifetimes are measured on a discrete scale; that is, the curing time of a particular disease measured in days, the survival time of cancer patients in months etc.

In the last decades, many authors developed discrete version of continuous distributions. The Normal distribution is one of the famous continuous distributions which Lisman and van Zuylen [3] and Kemp [4] studied its discrete version. Recently, Roy [5] presented another version of discrete normal distribution. Another versions of discrete distributions have been studied such as discrete exponential distribution (see Sato et al. [6]). In addition, Nekoukhou et al. [7] proposed a discrete analogue of the generalized exponential distribution which has been presented by Gupta and Kundu [8]. Weibull distribution tends to better represent life data and are commonly called lifetime distribution. Since, the continuous lifetime may often be counted as discrete random variable, Nakagawa and Osaki [9] presented a discrete Weibull distribution. Stein and Dattero [10] introduced a second type discrete Weibull distribution and a third one was developed by Padgett and Spurrier [11]. Recently, Jayakumar and Girish Babu [12] introduced a discrete analogue of Weibull geometric distribution. Also, they showed sub models of this distribution; that is, discrete Weibull, discrete Rayleigh and geometric distributions. A comprehensive discussion of many discrete versions of continuous distributions are found in Inusah and Kozubowski [13], Kozubowski and Inusah [14], Krishna and Pundir [15], Chakraborty et al. [16-18], Roy and Gupta [19], Roy[20] and Nekoukhou and Bidram [21].

One of the more applicable lifetime distributions is the Wiener process range distribution which has been obtained by Feller [22]. He used the method of images to derive the probability density function of this range (random variable). This range is the best to express about the change of price formula based on the assumption that stock price follow a Wiener process. It is given from the difference between the highest and the lowest stock price. Recently, Withers and Nadarajah [23] presented an expansion for this distribution by giving its cumulative distribution function and its quantiles. More recently, Teamah et al. [1] provided the truncated distribution of a Wiener process range and studied its various statistical properties.

In order to evaluate the behavior of some stocks over specified number of months (discrete points), the economists should know the difference between the highest and lowest price during each month (range in month ≡ point's value in this month). Generally, this difference is measured on a continuous scale over the time in each month, see Teamah et al. [1], Feller [22] and Withers and Nadarajah [23], but for the previous reason we propose to deal with a model has an integer values. Thus, in this paper we should find the discrete distribution of these point's value which consider as a Wiener process range in months. We will study some basic distributional properties including reliability properties, moments, stress-strength parameter and order statistics.

The paper is organized as follows. In Section 2, we introduce the Discrete Distribution for a Wiener Process Range (DDWPR). We study some basic distributional properties for DDWPR in Section 3. Also, we present the truncated version of DDWPR which is denoted by TDDWPR and its statistical properties in Section 4. Section 5 presented an application to a real data set. The ends of the paper with some concluding remarks and future works are discussed in Section 6.

## 2. Discrete Distribution for a Wiener Process Range (DDWPR)

In the time interval $(0,T)$ the range of the Wiener process $\{W(t); t \geq 0\}$ is $\bar{R}(T) = \sup_{(0,T)} W(t) - \inf_{(0,T)} W(t)$ and its probability density function (pdf) is given by:

$$\bar{f}_{\bar{R}(T)}(\bar{r}) = \left(\frac{2}{\pi}\right)^{\frac{1}{2}} \bar{r}^{-1}(2\pi)^{\frac{1}{2}}\left(\frac{\bar{r}T^{-\frac{1}{2}}}{2}\right)^{-1} \sum_{k=1}^{\infty} \exp\left[-\frac{(2k-1)^2 \pi^2}{8} \cdot \left(\frac{\bar{r}T^{-\frac{1}{2}}}{2}\right)^{-2}\right] \quad (1)$$

Feller [22]. Since, $\bar{R}(T)$ has a pdf (1), then we get:

$$P[\bar{R}(T) > h] = 1 - P[\bar{R}(T) \leq h] = 1 - \sum_{k=1}^{\infty}\left(\frac{8}{(2k-1)^2 \pi^2} + \frac{8T}{h^2}\right)\exp\left[-\frac{(2k-1)^2 \pi^2 T}{2h^2}\right]. \quad (2)$$

Set $R(T) = \lfloor \bar{R}(T) \rfloor$, where $\lfloor \ \rfloor$ denotes to the minimum integer value and let $r$ be a positive integer, then by definition of $R(T)$ we have $P[R(T) \geq r] = P[\lfloor \bar{R}(T) \rfloor \geq r]$. Since, $r$ is a an integer, then this is actually equal to $P[\bar{R}(T) \geq r]$ which equal to $P[R(T) \geq r]$. Consequently, from (2) we have

$$P[R(T) \geq r] = 1 - \sum_{k=1}^{\infty} \left( \frac{8}{(2k-1)^2 \pi^2} + \frac{8T}{r^2} \right) \exp\left[ -\frac{(2k-1)^2 \pi^2 T}{2r^2} \right]. \tag{3}$$

Since, the event $\{R(T) = r\}$ is the same as the event $\{R(T) \geq r\} \cap \{R(T) \geq r+1\}$, then the probability mass function (pmf) of DDWPR is given by,

$$f_{R(T)}(r) = P[R(T) = r] = P[R(T) \geq r] - P[R(T) \geq r+1]$$

$$= \sum_{k=1}^{\infty} \left[ \left( \frac{8}{(2k-1)^2 \pi^2} + \frac{8T}{(r+1)^2} \right) \exp\left[ -\frac{(2k-1)^2 \pi^2 T}{2(r+1)^2} \right] - \left( \frac{8}{(2k-1)^2 \pi^2} + \frac{8T}{r^2} \right) \exp\left[ -\frac{(2k-1)^2 \pi^2 T}{2r^2} \right] \right],$$

$$\tag{5}$$

$r = 0,1,2,\ldots$ and it is represented as in figure 1.

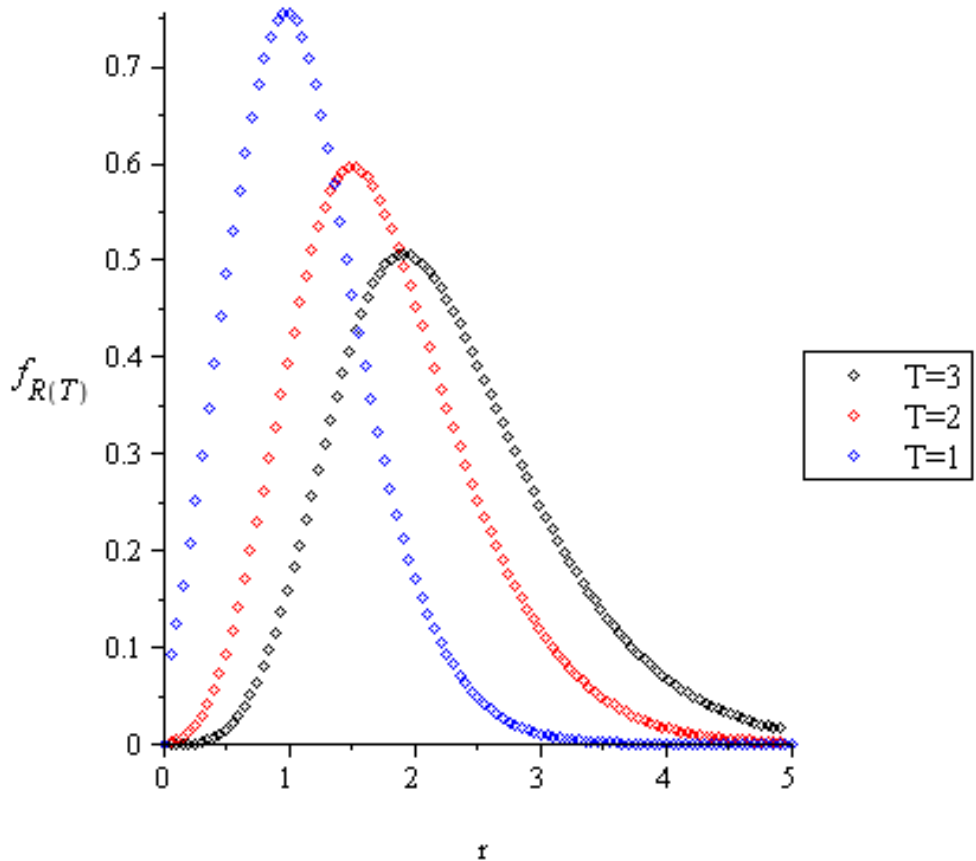

**Figure 1:** The probability mass function of DDWPR.

It is clear that :

$$\sum_{r=0}^{\infty} P[R(T)=r] = P[R(T) \geq 0] - P[R(T) \geq 1] + P[R(T) \geq 1] - P[R(T) \geq 2] + \ldots = P[R(T) \geq 0] = S_{R(T)}(0) = 1,$$

where $S_{R(T)}(.)$ is the survival function of the random variable $R(T)$.

## 2.1 Cumulative distribution function

The cumulative distribution function (cdf) of DDWPR is given by,

$$F_{R(T)}(r) = P[R(T) \leq r] = 1 - P[\overline{R}(T) \geq r] + P[R(T) = r] = 1 - P[R(T) \geq r+1]$$

$$= \sum_{k=1}^{\infty} \left[ \left( \frac{8}{(2k-1)^2 \pi^2} + \frac{8T}{(r+1)^2} \right) \exp\left[ -\frac{(2k-1)^2 \pi^2 T}{2(r+1)^2} \right] \right], \quad (6)$$

$r = 0,1,2,\ldots$ and it is represented as in figure 2.

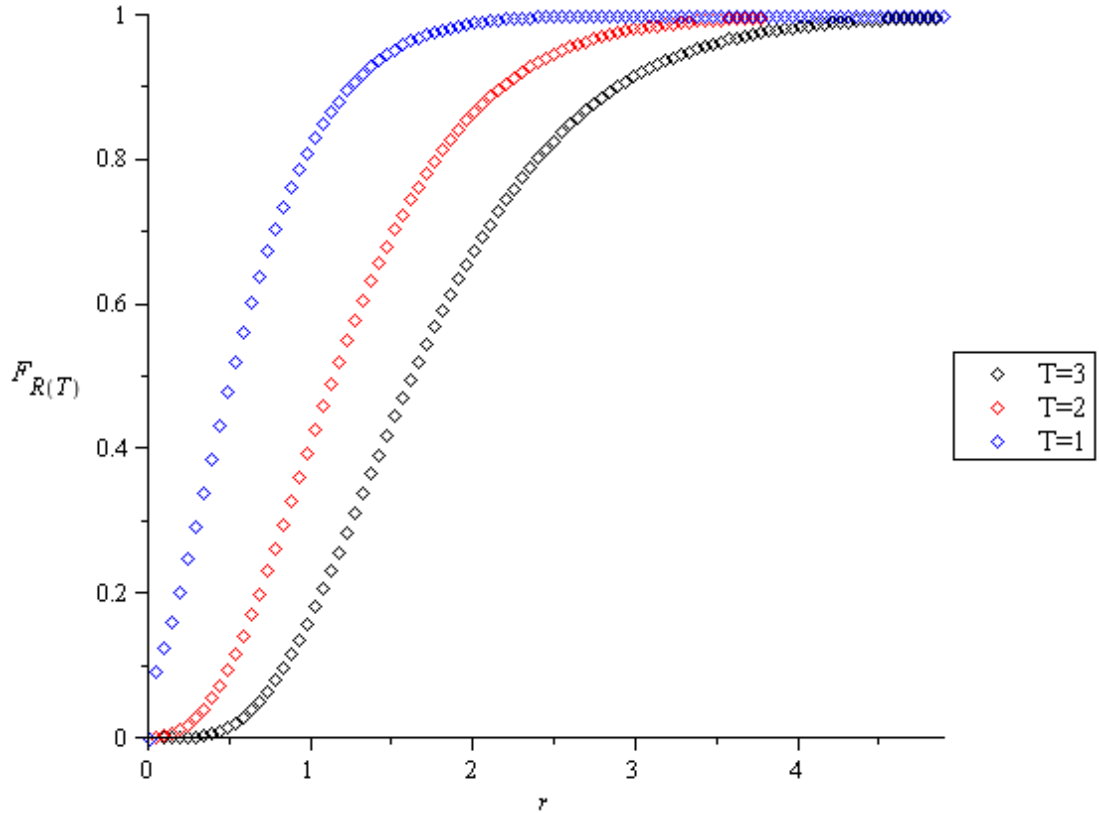

**Figure 2**: Cumulative distribution function of DDWPR.

We note that $F_{R(T)}(\infty) = 1$ where $\sum_{k=1}^{\infty} \frac{8}{(2k-1)^2 \pi^2} = 1$.

## 2.2 Survival function

In discrete distributions, it is easy to show that the survival function $S_{R(T)}(.) = 1 - F_{R(T)}(r)$ is decreasing by increasing the value of $T$, see figure 3.

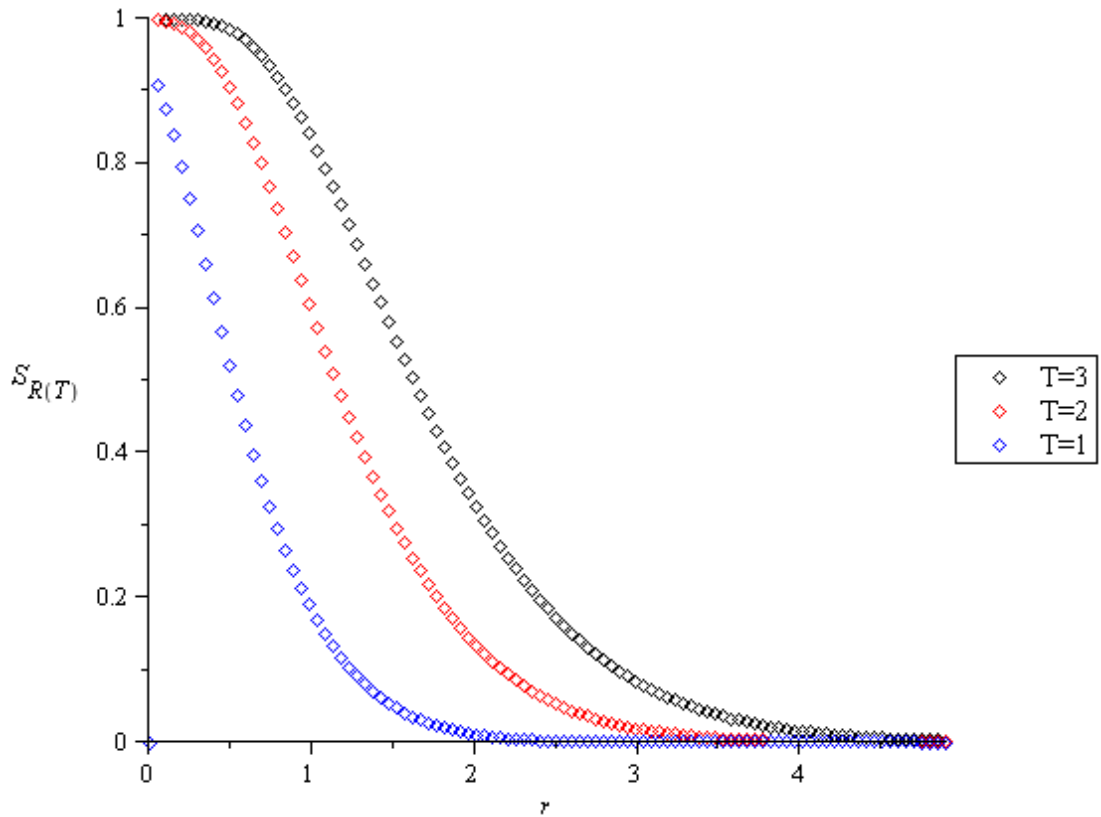

**Figure 3:** Survival function of DDWPR.

## 3. Some basic distributional properties of DDWPR

In our problem we need to find some basic measures for stock price over specified periods of time (days or months or years). Thus, in this section we study some various statistical properties of the DDWPR including reliability properties, moments, stress-strength parameter, Bonferroni curve, Lorenz curve and Gini's index.

### 3.1 Reliability properties of DDWPR

In the economics science the risk rate (hazard rate) is influenced by the swings between fall and rise much of the stock price during the time period $(0,T)$. If we divide this period into specified periods of time (days or months or years), then we want to find the discrete hazard rate which is more applicable in reliability theory. Thus, as in Shaked et al. [24] the hazard rate function of DDWPR is given by,

$$h_{R(T)}(r) = \frac{P[R(T) = r]}{P[R(T) \geq r]}$$

$$= \frac{\sum_{k=1}^{\infty}\left[\left(\frac{8}{(2k-1)^2 \pi^2} + \frac{8T}{(r+1)^2}\right)\exp\left[-\frac{(2k-1)^2 \pi^2 T}{2(r+1)^2}\right] - \left(\frac{8}{(2k-1)^2 \pi^2} + \frac{8T}{r^2}\right)\exp\left[-\frac{(2k-1)^2 \pi^2 T}{2r^2}\right]\right]}{1 - \sum_{k=1}^{\infty}\left[\left(\frac{8}{(2k-1)^2 \pi^2} + \frac{8T}{(r+1)^2}\right)\exp\left[-\frac{(2k-1)^2 \pi^2 T}{2(r+1)^2}\right]\right]},$$

(7)

where $P[R(T) \geq r] > 0$ and it is represented as in figure 4.

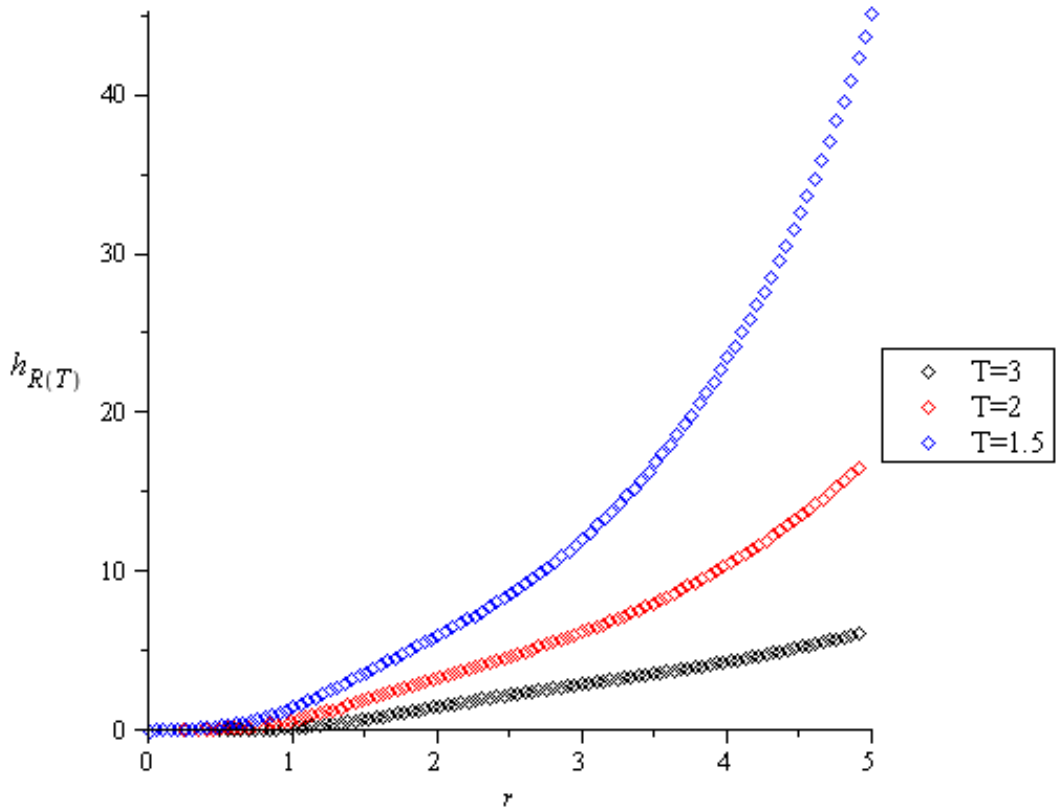

**Figure 4:** Hazard rate function of DDWPR.

It is clear that, at $r = 0$ we found:

$$h_{R(T)}(0) = f_{R(T)}(0) = \sum_{k=1}^{\infty}\left(\frac{8}{(2k-1)^2 \pi^2} + 8T\right)\exp\left[-\frac{(2k-1)^2 \pi^2 T}{2}\right]. \quad (8)$$

The mean residual life function (MRL) of DDWPR is given by,

$$L_{R(T)}(r) = E[(R(T)-r)|R(T) \geq r] = \frac{\sum_{j \geq r} j f_{R(T)}(j)}{\sum_{j \geq r} f_{R(T)}(j)} - r = \frac{\sum_{j \geq r} S_{R(T)}(j)}{S_{R(T)}(r)} = \sum_{j \geq r} \prod_{i=r}^{j}(1-h_{R(T)}(i))$$

$$= \sum_{j \geq r} \prod_{i=r}^{j} \left[ \frac{1+\sum_{k=1}^{\infty}\left[\left(\frac{8}{(2k-1)^2\pi^2}+\frac{8T}{i^2}\right)\exp\left[-\frac{(2k-1)^2\pi^2 T}{2i^2}\right] - 2\left(\frac{8}{(2k-1)^2\pi^2}+\frac{8T}{(i+1)^2}\right)\exp\left[-\frac{(2k-1)^2\pi^2 T}{2(i+1)^2}\right]\right]}{1-\sum_{k=1}^{\infty}\left[\left(\frac{8}{(2k-1)^2\pi^2}+\frac{8T}{(i+1)^2}\right)\exp\left[-\frac{(2k-1)^2\pi^2 T}{2(i+1)^2}\right]\right]} \right],$$

(9)

where $r = 0,1,2,\ldots$ . On the other hand, Roy and Gupta [19] presented another formula for MRL,

$$\mu_{R(T)}(r) = E[(R(T)-r)|R(T) \geq r] = L_{R(T)}(r+1) + 1$$

$$= 1 + \sum_{j \geq r+1} \prod_{i=r+1}^{j} \left[ \frac{1+\sum_{k=1}^{\infty}\left[\left(\frac{8}{(2k-1)^2\pi^2}+\frac{8T}{i^2}\right)\exp\left[-\frac{(2k-1)^2\pi^2 T}{2i^2}\right] - 2\left(\frac{8}{(2k-1)^2\pi^2}+\frac{8T}{(i+1)^2}\right)\exp\left[-\frac{(2k-1)^2\pi^2 T}{2(i+1)^2}\right]\right]}{1-\sum_{k=1}^{\infty}\left[\left(\frac{8}{(2k-1)^2\pi^2}+\frac{8T}{(i+1)^2}\right)\exp\left[-\frac{(2k-1)^2\pi^2 T}{2(i+1)^2}\right]\right]} \right],$$

(10)

where $r = 0,1,2,\ldots$. If we consider $r = 0$, then the MRL function is equal to the mean of the lifetime distribution (see Jayakumar and Girish Babu [12]) where $L_{R(T)}(0) = \mu$. Thus, we have,

$$\mu_{R(T)}(0) = \frac{\mu}{1-f_{R(T)}(0)} = \frac{\mu}{1-\sum_{k=1}^{\infty}\left(\frac{8}{(2k-1)^2\pi^2}+8T\right)\exp\left[-\frac{(2k-1)^2\pi^2 T}{2}\right]}.$$

(11)

Also, the reversed hazard rate function of DDWPR is given by,

$$h^*_{R(T)}(r) = P[(R(T)=r)/(R(T) \leq r)] = \frac{P[R(T)=r]}{P[R(T) \leq r]}$$

$$= \frac{\sum_{k=1}^{\infty}\left[\left(\frac{8}{(2k-1)^2\pi^2}+\frac{8T}{(r+1)^2}\right)\exp\left[-\frac{(2k-1)^2\pi^2 T}{2(r+1)^2}\right]-\left(\frac{8}{(2k-1)^2\pi^2}+\frac{8T}{r^2}\right)\exp\left[-\frac{(2k-1)^2\pi^2 T}{2r^2}\right]\right]}{\sum_{k=1}^{\infty}\left[\left(\frac{8}{(2k-1)^2\pi^2}+\frac{8T}{(r+1)^2}\right)\exp\left[-\frac{(2k-1)^2\pi^2 T}{2(r+1)^2}\right]\right]},$$

(12)

see figure 5.

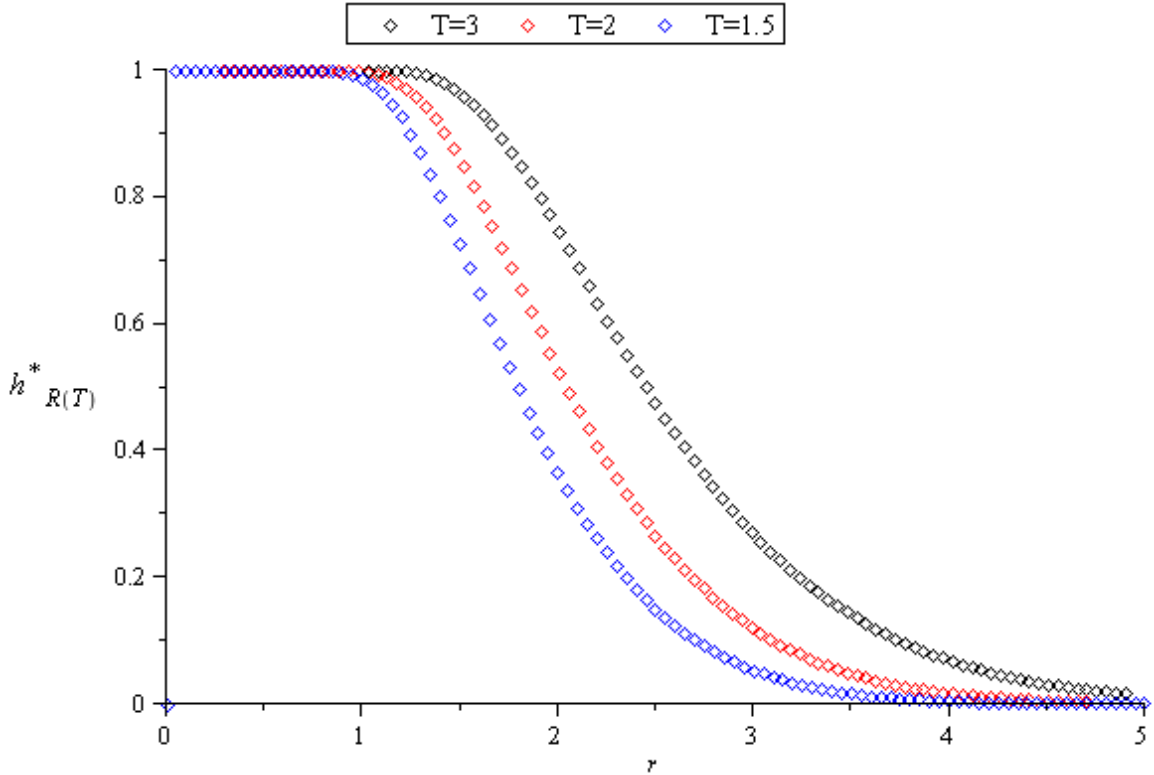

**Figure 5:** Reversed hazard rate function of DDWPR.

Also, the second rate of failure of DDWPR is given by, $h_{R(T)}^{**}(r)=\log\left\{\frac{S_{R(T)}(r)]}{S_{R(T)}(r+1)}\right\}$, (see Jayakumar and Girish Babu [12]) and it is represented as in figure 6.

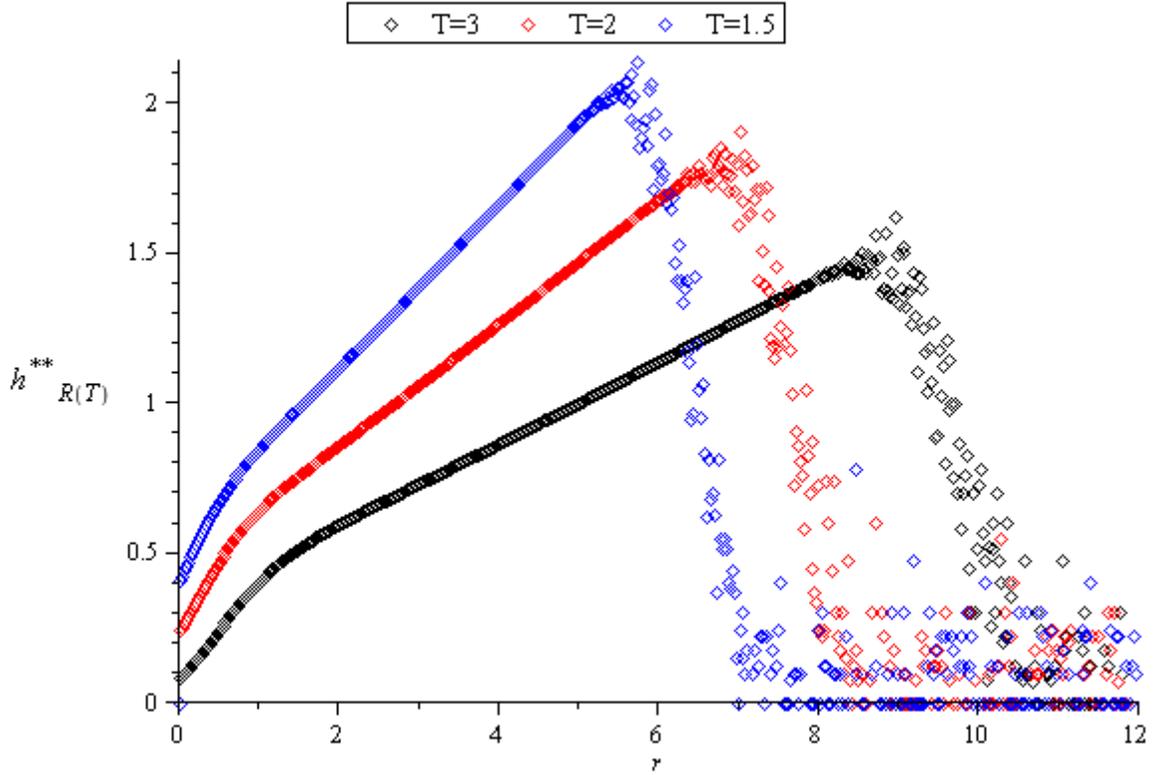

**Figure 6:** Second rate function of DDWPR.

## 3.2 DDWPR quantiles's and its random number generation

Wiener range distribution is one of the important real lifetime distributions which has no parameters, because this distribution is derived from an underlying process which has no parameters (standard Gaussian). Sometimes the economists need to conduct a deeper study on the stock price behavior at specific points. They divide the domain of this distribution into equal probability parts or the statistical sample into equal size groups by using quantile. As in Rohatgi and Saleh [25], we let the point $r_u$ is known as the $u^{th}$ quantile of $R$ if it satisfies, $P[R(T) \leq r_u] \geq u$ and $P[R(T) \geq r_u] \geq 1-u$.

**Remark 1.** *If $0 < a_i < 1$, $i = 1,2,...,n$ then $\prod_{i=1}^{n} a_i < \sum_{i=1}^{n} a_i$.*

**Theorem 1.** *For the DDWPR, the $u^{th}$ quantile $\psi(u)$ is given by, $\psi(u) = \lceil r_u \rceil$ where $\lceil r_u \rceil \geq r_u$ and $r_u$ is given from solving the following equation:*

$$P[R(T) \geq r_u] + \left(\frac{8T}{(r_u+1)^2} + 1\right)\left(\frac{1}{2}\exp\left[\frac{3\pi^2 T}{2(r_u+1)^2}\right]\right)\left(\exp\left[-\frac{2\pi^2 T}{(r_u+1)^2}\right]\right)^{\frac{3}{4}} \vartheta\left[2, 0, \exp\left[-\frac{2\pi^2 T}{(r_u+1)^2}\right]\right] = 0. \quad (13)$$

**Proof.** Since, $P[R(T) \leq r_u] \geq u$ then we have,

$$\sum_{k=1}^{\infty}\left[\left(\frac{8}{(2k-1)^2 \pi^2} + \frac{8T}{(r_u+1)^2}\right)\exp\left[-\frac{(2k-1)^2 \pi^2 T}{2(r_u+1)^2}\right]\right] \geq u,$$

$$\Rightarrow \sum_{k=1}^{\infty}\left[\left(\frac{8}{(2k-1)^2 \pi^2}\right)\exp\left[-\frac{(2k-1)^2 \pi^2 T}{2(r_u+1)^2}\right]\right] + \sum_{k=1}^{\infty}\left[\left(\frac{8T}{(r_u+1)^2}\right)\exp\left[-\frac{(2k-1)^2 \pi^2 T}{2(r_u+1)^2}\right]\right] \geq u.$$

For any value of $k = 1, 2, \ldots, \infty$, we have $\frac{8}{(2k-1)^2 \pi^2} < 1$ and also $\exp\left[-\frac{(2k-1)^2 \pi^2 T}{2(r_u+1)^2}\right] < 1$ where $r_u$ and $T$ are fixed values (discrete case). Consequently, by using remark 1 we get,

$$\sum_{k=1}^{\infty}\left(\frac{8}{(2k-1)^2 \pi^2}\right) + \sum_{k=1}^{\infty}\exp\left[-\frac{(2k-1)^2 \pi^2 T}{2(r_u+1)^2}\right] + \sum_{k=1}^{\infty}\left[\left(\frac{8T}{(r_u+1)^2}\right)\exp\left[-\frac{(2k-1)^2 \pi^2 T}{2(r_u+1)^2}\right]\right]$$

$$< \sum_{k=1}^{\infty}\left[\left(\frac{8}{(2k-1)^2 \pi^2}\right)\exp\left[-\frac{(2k-1)^2 \pi^2 T}{2(r_u+1)^2}\right]\right] + \sum_{k=1}^{\infty}\left[\left(\frac{8T}{(r_u+1)^2}\right)\exp\left[-\frac{(2k-1)^2 \pi^2 T}{2(r_u+1)^2}\right]\right]$$

$$\geq u,$$

leads to,

$$\sum_{k=1}^{\infty}\left(\frac{8}{(2k-1)^2 \pi^2}\right) + \sum_{k=1}^{\infty}\exp\left[-\frac{(2k-1)^2 \pi^2 T}{2(r_u+1)^2}\right] + \sum_{k=1}^{\infty}\left[\left(\frac{8T}{(r_u+1)^2}\right)\exp\left[-\frac{(2k-1)^2 \pi^2 T}{2(r_u+1)^2}\right]\right] \leq u$$

$$\Rightarrow 1 + \left(\frac{8T}{(r_u+1)^2} + 1\right)\left(\frac{1}{2}\exp\left[\frac{3\pi^2 T}{2(r_u+1)^2}\right]\right)\left(\exp\left[-\frac{2\pi^2 T}{(r_u+1)^2}\right]\right)^{\frac{3}{4}} \vartheta\left[2, 0, \exp\left[-\frac{2\pi^2 T}{(r_u+1)^2}\right]\right] \leq u,$$

where $\vartheta$ is *EllipticTheta* function.

Also, from $P[R(T) \geq r_u] \geq 1-u$. Thus, using (3) the integer value of $r_u$ is given from solving the following equation:

$$P[R(T) \geq \lceil r_u \rceil] + \left(\frac{8T}{(\lceil r_u \rceil+1)^2} + 1\right)\left(\frac{1}{2}\exp\left[\frac{3\pi^2 T}{2(\lceil r_u \rceil+1)^2}\right]\right)\left(\exp\left[-\frac{2\pi^2 T}{(\lceil r_u \rceil+1)^2}\right]\right)^{\frac{3}{4}} \vartheta\left[2, 0, \exp\left[-\frac{2\pi^2 T}{(\lceil r_u \rceil+1)^2}\right]\right] = 0.$$

This completes the proof. □

Due to the random changes in stock prices during intermittent period, the random numbers (integer) are important for simulating these events and other quantities. This application may require uniformly distributed numbers. Thus, if the random number $U$ is drawn from a uniform distribution on (0, 1), then a random number $R(T)$ following DDWPR is given from (13). In particular, the median is given solving the following equation:

$$P[R(T) \geq \lceil r_{0.5} \rceil] + \left(\frac{8T}{(\lceil r_{0.5} \rceil+1)^2}+1\right)\left(\frac{1}{2}\exp\left[\frac{3\pi^2 T}{2(\lceil r_{0.5} \rceil+1)^2}\right]\right)\left(\exp\left[-\frac{2\pi^2 T}{(\lceil r_{0.5} \rceil+1)^2}\right]\right)^{\frac{3}{4}} \cdot \mathcal{G}\left[2,0,\exp\left[-\frac{2\pi^2 T}{(\lceil r_{0.5} \rceil+1)^2}\right]\right] = 0.$$

(14)

### 3.3 Moments of DDWPR

The $q^{th}$ moment about origin is given by,

$$\mu'_q = \sum_{q=0}^{\infty}\left[r^q \sum_{k=1}^{\infty}\left[\left(\frac{8}{(2k-1)^2 \pi^2}+\frac{8T}{(r+1)^2}\right)\exp\left[-\frac{(2k-1)^2 \pi^2 T}{2(r+1)^2}\right]-\left(\frac{8}{(2k-1)^2 \pi^2}+\frac{8T}{r^2}\right)\exp\left[-\frac{(2k-1)^2 \pi^2 T}{2r^2}\right]\right]\right].$$

It is difficult to calculate the above series, but one can computes it numerically. Thus, we can obtain the values of moments for given values of $T$ by using Maple programming as in the following Table 1.

**Table 1.** Moments, skewness and kurtosis for various values of $T$.

| $T$ | Raw moments | Central moments | Skewness | Kurtosis |
|---|---|---|---|---|
| 1 | $\mu'_1 = 1.129581778$<br>$\mu'_2 = 1.536671058$<br>$\mu'_3 = 2.4247010$<br>$\mu'_4 = 4.6151$ | $\mu_2 = 0.260716259$<br>$\mu_3 = 25.72779723$<br>$\mu_4 = 0.539663483$ | 193.2637627 | 4.939383864 |
| 2 | $\mu'_1 = 1.747634563$<br>$\mu'_2 = 3.598024786$<br>$\mu'_3 = 8.512482209$<br>$\mu'_4 = 22.80606816$ | $\mu_2 = 0.543798220$<br>$\mu_3 = 1223.317372$<br>$\mu_4 = 1.24943324$ | 3050.580364 | 1.225104865 |
| 3 | $\mu'_1 = 2.263047990$<br>$\mu'_2 = 5.890473359$<br>$\mu'_3 = 17.32622321$<br>$\mu'_4 = 57.0669$ | $\mu_2 = 0.769087154$<br>$\mu_3 = 10379.93043$<br>$\mu_4 = 2.54514560$ | 15389.73376 | 1.302902616 |

## 3.4. Order Statistics of DDWPR

Order statistics among the most fundamental tools in non-parametric statistics and inference have a remarkable importance in economic due to its great applicability. It is an important methodology to find the exact distribution of the maximum and the minimum of the prices-path.

Let $R_{1:n} \leq R_{2:n} \leq ... \leq R_{n:n}$ denote the order statistics of a random sample $R_1, R_2, ..., R_n$ from the DDWPR, then the pmf of the p*th* order statistic $R_{p:n}$ is,

$$f_{R(T)(p:n)} = \frac{n!}{(p-1)!(n-p)!} \left(F_{R(T)}(r)\right)^{p-1} \left(1-F_{R(T)}(r)\right)^{n-p} f_{R(T)}(r),$$

$$= \frac{n!}{(p-1)!(n-p)!} \left( \sum_{k=1}^{\infty} \left[ \left( \frac{8}{(2k-1)^2 \pi^2} + \frac{8T}{(r+1)^2} \right) \exp\left[ -\frac{(2k-1)^2 \pi^2 T}{2(r+1)^2} \right] \right] \right)^{p-1}$$

$$\times \left( 1 - \sum_{k=1}^{\infty} \left[ \left( \frac{8}{(2k-1)^2 \pi^2} + \frac{8T}{(r+1)^2} \right) \exp\left[ -\frac{(2k-1)^2 \pi^2 T}{2(r+1)^2} \right] \right] \right)^{n-p}$$

$$\times \left( \sum_{k=1}^{\infty} \left[ \left( \frac{8}{(2k-1)^2 \pi^2} + \frac{8T}{(r+1)^2} \right) \exp\left[ -\frac{(2k-1)^2 \pi^2 T}{2(r+1)^2} \right] - \left( \frac{8}{(2k-1)^2 \pi^2} + \frac{8T}{r^2} \right) \exp\left[ -\frac{(2k-1)^2 \pi^2 T}{2r^2} \right] \right] \right),$$

(15)

Also, the distribution function of $R_{p:n}$ is,

$$F_{R(T)(p:n)}(r) = \sum_{i=p}^{n} \binom{n}{i} \left(F_{R(T)}(r)\right)^i \left(1-F_{R(T)}(r)\right)^{n-i}$$

$$= \sum_{i=p}^{n} \binom{n}{i} \left( \sum_{k=1}^{\infty} \left[ \left( \frac{8}{(2k-1)^2 \pi^2} + \frac{8T}{(r+1)^2} \right) \exp\left[ -\frac{(2k-1)^2 \pi^2 T}{2(r+1)^2} \right] \right] \right)^i$$

$$\times \left( 1 - \sum_{k=1}^{\infty} \left[ \left( \frac{8}{(2k-1)^2 \pi^2} + \frac{8T}{(r+1)^2} \right) \exp\left[ -\frac{(2k-1)^2 \pi^2 T}{2(r+1)^2} \right] \right] \right)^{n-i}.$$

(16)

In addition, we can get the pmf of the minimum difference between the prices at knowing time $T$ as:

$$f_{R(T)(1:n)} = n\left(1 - \sum_{k=1}^{\infty}\left[\left(\frac{8}{(2k-1)^2\pi^2} + \frac{8T}{(r+1)^2}\right)\exp\left[-\frac{(2k-1)^2\pi^2T}{2(r+1)^2}\right]\right]\right)^{n-1}$$

$$\times \left(\sum_{k=1}^{\infty}\left[\left(\frac{8}{(2k-1)^2\pi^2} + \frac{8T}{(r+1)^2}\right)\exp\left[-\frac{(2k-1)^2\pi^2T}{2(r+1)^2}\right] - \left(\frac{8}{(2k-1)^2\pi^2} + \frac{8T}{r^2}\right)\exp\left[-\frac{(2k-1)^2\pi^2T}{2r^2}\right]\right]\right).$$

(17)

and the pmf of the maximum difference is,

$$f_{R(T)(n:n)} = n\left(\sum_{k=1}^{\infty}\left[\left(\frac{8}{(2k-1)^2\pi^2} + \frac{8T}{(r+1)^2}\right)\exp\left[-\frac{(2k-1)^2\pi^2T}{2(r+1)^2}\right]\right]\right)^{n-1}$$

$$\times \left(\sum_{k=1}^{\infty}\left[\left(\frac{8}{(2k-1)^2\pi^2} + \frac{8T}{(r+1)^2}\right)\exp\left[-\frac{(2k-1)^2\pi^2T}{2(r+1)^2}\right] - \left(\frac{8}{(2k-1)^2\pi^2} + \frac{8T}{r^2}\right)\exp\left[-\frac{(2k-1)^2\pi^2T}{2r^2}\right]\right]\right),$$

(18)

### 3.5. Stress-strength parameter of DDWPR

In discrete case, if the random variable $R(T)$ which describes the changing of stock price is the strength of a component which is subjected to a random stress $Z$, then the stress-strength model is defined as,

$$P(R(T) > Z) = \sum_{r=0}^{\infty} f_{R(T)}(r) F_{R(T)}(r)$$

If $R(T)$ has $DDWPR(T_1)$ and $Z$ has $DDWPR(T_2)$ and using (5), (6), then we have

$$\Delta = \sum_{r=0}^{\infty}\sum_{k=1}^{\infty}\left[\left(\frac{8}{(2k-1)^2\pi^2} + \frac{8T}{(r+1)^2}\right)\exp\left[-\frac{(2k-1)^2\pi^2T}{2(r+1)^2}\right] - \left(\frac{8}{(2k-1)^2\pi^2} + \frac{8T}{r^2}\right)\exp\left[-\frac{(2k-1)^2\pi^2T}{2r^2}\right]\right]$$

$$\times \sum_{k=1}^{\infty}\left[\left(\frac{8}{(2k-1)^2\pi^2} + \frac{8T}{(r+1)^2}\right)\exp\left[-\frac{(2k-1)^2\pi^2T}{2(r+1)^2}\right]\right],$$

where $\Delta = P(R(T) > Z)$ is stress-strength parameter.

### 4. Truncated Discrete Distribution for a Wiener Process Range (TDDWPR)

In some real life problems we need to delete some values of $R(T)$ then we use the truncation method. The main idea of the truncation method is to distribute the probability of the deleted part on the other part. Thus, this method is an important methodology in different fields of sciences, in particular communication networks and finance. etc. There are three kinds of truncation, right, left and double truncation. In an earlier work, this method has been studied extensively for many lifetime distributions by some authors, for example Zaninetti [26], Ali and Nadarajah [27], Nadarajah [28], Pender [29] and Chattopadhyay et al. [30]. Recently, Teamah et al. [1] studied the truncated distribution of a Wiener process range and studied its various statistical properties.

As in Teamah et al. [1], if $R(T)$ is a random variable with probability mass function (5), define $\tilde{R}(T)$ as a corresponding double truncated (truncation from left and right) of $R(T)$ with the probability mass function $q_{\tilde{R}(T)}(r)$ is given by: Then, the pmf of double truncated of $R(T)$ is given by:

$$q_{\tilde{R}(T)}(r) = \frac{f_{R(T)}(r)}{F_{R(T)}(b) - F_{R(T)}(a)}$$

$$= \frac{\sum_{k=1}^{\infty}\left[\left(\frac{8}{(2k-1)^2\pi^2} + \frac{8T}{(r+1)^2}\right)\exp\left[-\frac{(2k-1)^2\pi^2 T}{2(r+1)^2}\right] - \left(\frac{8}{(2k-1)^2\pi^2} + \frac{8T}{r^2}\right)\exp\left[-\frac{(2k-1)^2\pi^2 T}{2r^2}\right]\right]}{\sum_{k=1}^{\infty}\left[\left(\frac{8}{(2k-1)^2\pi^2} + \frac{8T}{(b+1)^2}\right)\exp\left[-\frac{(2k-1)^2\pi^2 T}{2(b+1)^2}\right]\right] - \sum_{k=1}^{\infty}\left[\left(\frac{8}{(2k-1)^2\pi^2} + \frac{8T}{(a+1)^2}\right)\exp\left[-\frac{(2k-1)^2\pi^2 T}{2(a+1)^2}\right]\right]},$$

(19)

$r = a,\ldots,a+i,\ldots,b,\ i$ is a non negative real number. Figure 7 represents the pmf of double truncated of $R(T)$. In addition, the cdf of TDDWPR is given by:

$$Q_{\tilde{R}(T)}(r) = \frac{\sum_{k=1}^{\infty}\left[\left(\frac{8}{(2k-1)^2\pi^2} + \frac{8T}{(r+1)^2}\right)\exp\left[-\frac{(2k-1)^2\pi^2 T}{2(r+1)^2}\right] - \left(\frac{8}{(2k-1)^2\pi^2} + \frac{8T}{(a+1)^2}\right)\exp\left[-\frac{(2k-1)^2\pi^2 T}{2(a+1)^2}\right]\right]}{\sum_{k=1}^{\infty}\left[\left(\frac{8}{(2k-1)^2\pi^2} + \frac{8T}{(b+1)^2}\right)\exp\left[-\frac{(2k-1)^2\pi^2 T}{2(b+1)^2}\right] - \left(\frac{8}{(2k-1)^2\pi^2} + \frac{8T}{(a+1)^2}\right)\exp\left[-\frac{(2k-1)^2\pi^2 T}{2(a+1)^2}\right]\right]},$$

(20)

see figure 8. Also, the survival function of TDDWPR is $\Omega_{\tilde{R}(T)}(r) = 1 - G_{\tilde{R}(T)}(r)$ and it is decreasing by increasing the value of $T$, see figure 9.

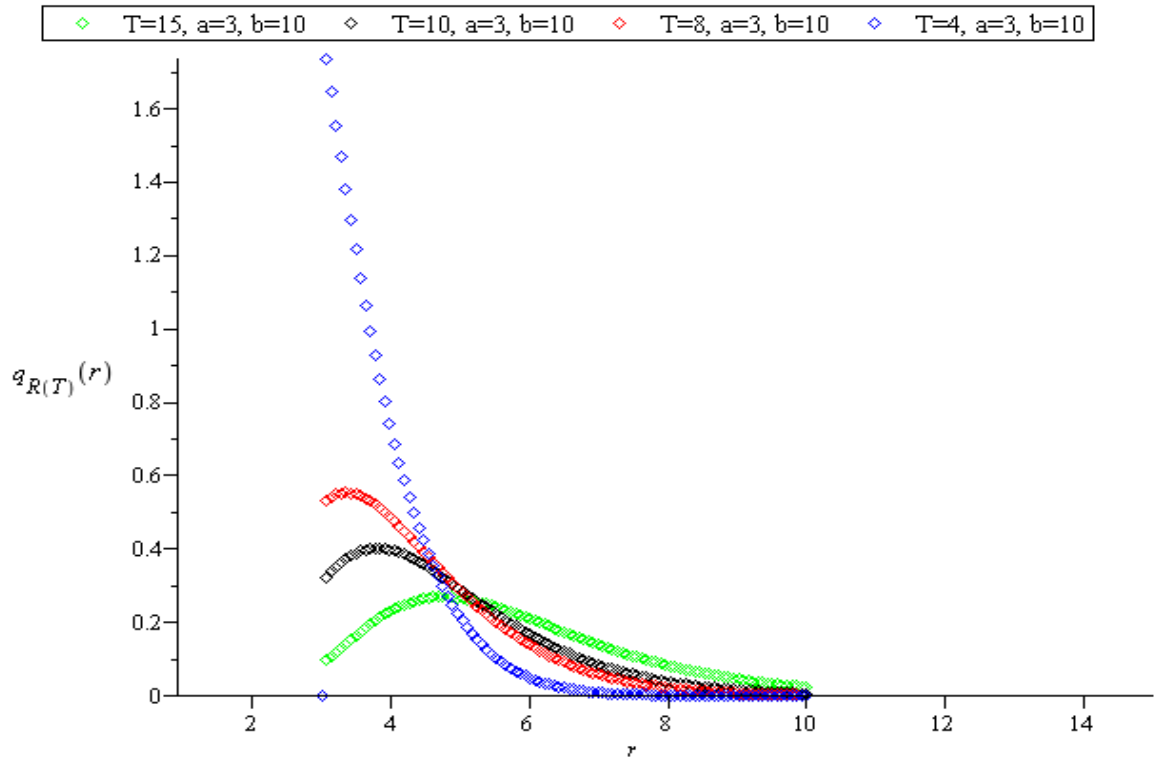

**Figure 7:** The probability mass function of TDDWPR.

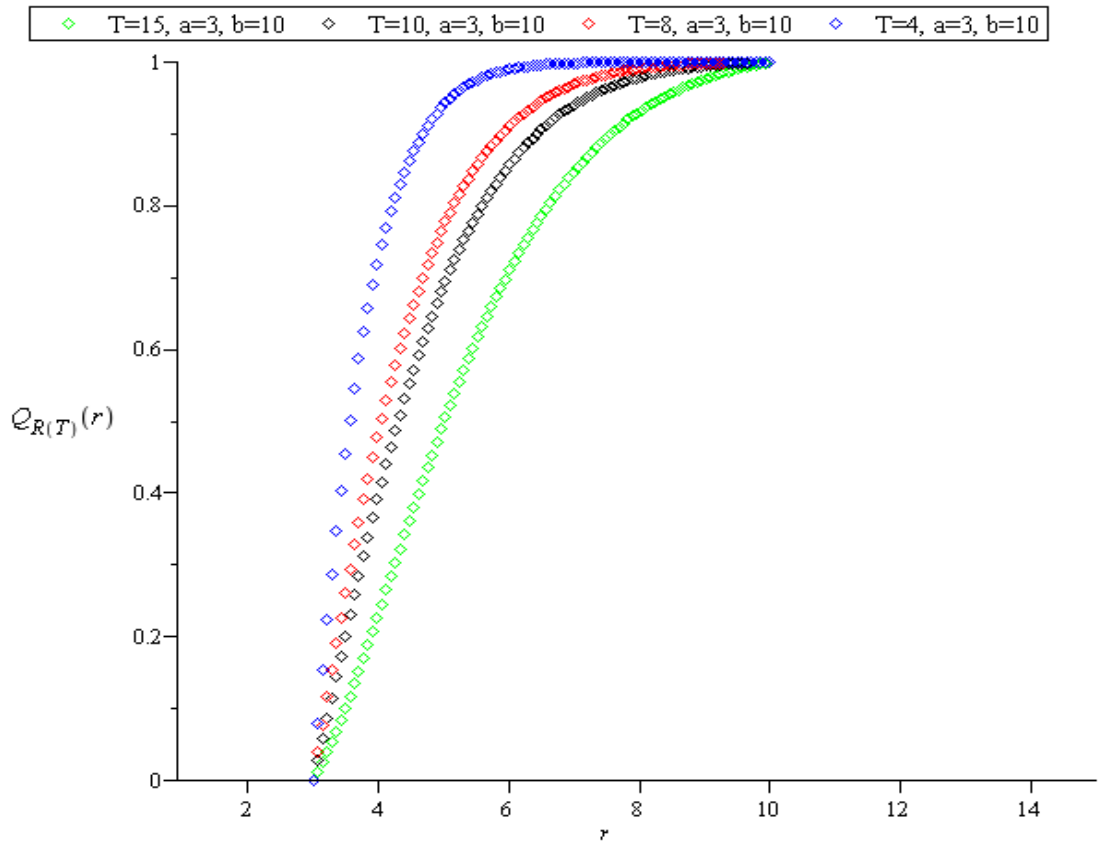

**Figure 8:** The cumulative distribution function of TDDWPR.

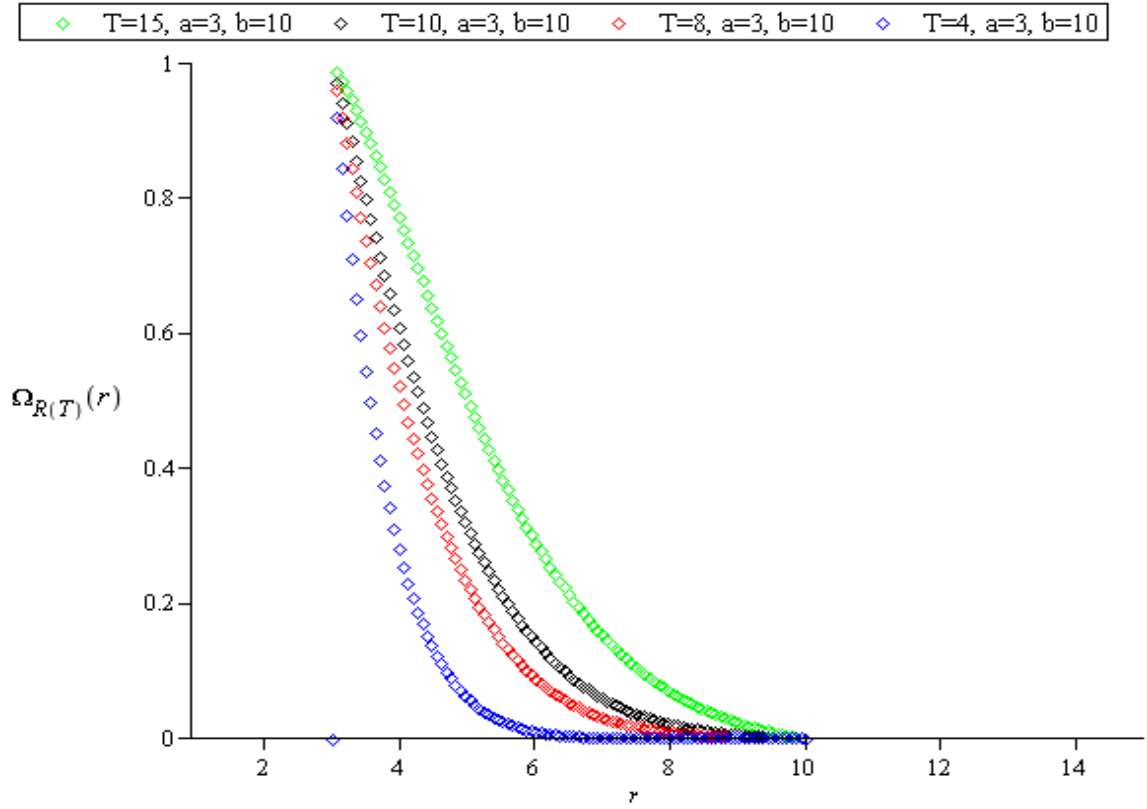

**Figure 9:** The survival function of TDDWPR.

## 4.1 Reliability properties of TDDWPR

Since, the risk rate (hazard rate) is an important measure for the stock price and it is influenced by the swings between fall and rise much of the stock price during the time period $(0,T)$. We get the hazard rate function of TDDWPR as follows:

$$\varphi_{\widetilde{R}(T)}(r) = q_{\widetilde{R}(T)}\{1 - Q_{\widetilde{R}(t)}(r)\}^{-1}$$

$$= \frac{\sum_{k=1}^{\infty}\left[\left(\frac{8}{(2k-1)^2\pi^2} + \frac{8T}{(r+1)^2}\right)\exp\left[-\frac{(2k-1)^2\pi^2 T}{2(r+1)^2}\right] - \left(\frac{8}{(2k-1)^2\pi^2} + \frac{8T}{r^2}\right)\exp\left[-\frac{(2k-1)^2\pi^2 T}{2r^2}\right]\right]}{\sum_{k=1}^{\infty}\left[\left(\frac{8}{(2k-1)^2\pi^2} + \frac{8T}{(b+1)^2}\right)\exp\left[-\frac{(2k-1)^2\pi^2 T}{2(b+1)^2}\right]\right] - \sum_{k=1}^{\infty}\left[\left(\frac{8}{(2k-1)^2\pi^2} + \frac{8T}{(a+1)^2}\right)\exp\left[-\frac{(2k-1)^2\pi^2 T}{2(a+1)^2}\right]\right]}$$
$$\div \left\{1 - \frac{\sum_{k=1}^{\infty}\left[\left(\frac{8}{(2k-1)^2\pi^2} + \frac{8T}{(r+1)^2}\right)\exp\left[-\frac{(2k-1)^2\pi^2 T}{2(r+1)^2}\right] - \left(\frac{8}{(2k-1)^2\pi^2} + \frac{8T}{(a+1)^2}\right)\exp\left[-\frac{(2k-1)^2\pi^2 T}{2(a+1)^2}\right]\right]}{\sum_{k=1}^{\infty}\left[\left(\frac{8}{(2k-1)^2\pi^2} + \frac{8T}{(b+1)^2}\right)\exp\left[-\frac{(2k-1)^2\pi^2 T}{2(b+1)^2}\right] - \left(\frac{8}{(2k-1)^2\pi^2} + \frac{8T}{(a+1)^2}\right)\exp\left[-\frac{(2k-1)^2\pi^2 T}{2(a+1)^2}\right]\right]}\right\}.$$

(21)

and it is represented as in figure 10. Also, figure 11 shows the reversed hazard rate function which is given by:

$$\varphi^*_{\tilde{R}(T)}(r) = q_{\tilde{R}(T)}\{Q_{\tilde{R}(t)}(r)\}^{-1}$$

$$= \frac{\sum_{k=1}^{\infty}\left[\left(\frac{8}{(2k-1)^2\pi^2} + \frac{8T}{(r+1)^2}\right)\exp\left[-\frac{(2k-1)^2\pi^2 T}{2(r+1)^2}\right] - \left(\frac{8}{(2k-1)^2\pi^2} + \frac{8T}{r^2}\right)\exp\left[-\frac{(2k-1)^2\pi^2 T}{2r^2}\right]\right]}{\sum_{k=1}^{\infty}\left[\left(\frac{8}{(2k-1)^2\pi^2} + \frac{8T}{(r+1)^2}\right)\exp\left[-\frac{(2k-1)^2\pi^2 T}{2(r+1)^2}\right] - \left(\frac{8}{(2k-1)^2\pi^2} + \frac{8T}{(a+1)^2}\right)\exp\left[-\frac{(2k-1)^2\pi^2 T}{2(a+1)^2}\right]\right]}. $$

(22)

In addition, the second rate function $\varphi^{**}_{\tilde{R}(T)}(r)$ is given by: $\varphi^{**}_{\tilde{R}(T)}(r) = \log\left\{\frac{\Omega_{\tilde{R}(T)}(r)}{\Omega_{\tilde{R}(T)}(r+1)}\right\}$ and it is represented as in figure 12.

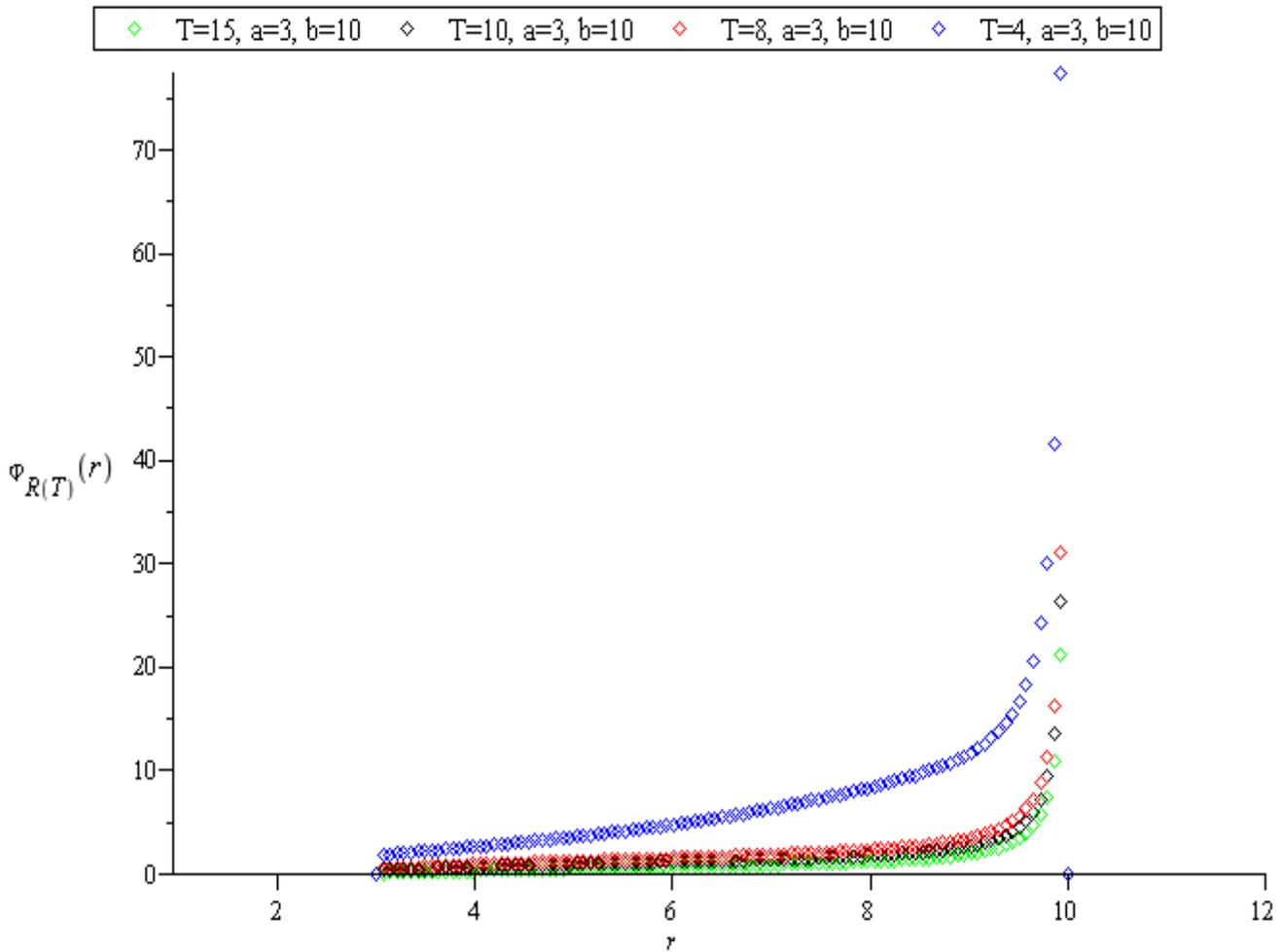

**Figure 10:** The hazard rate function of TDDWPR.

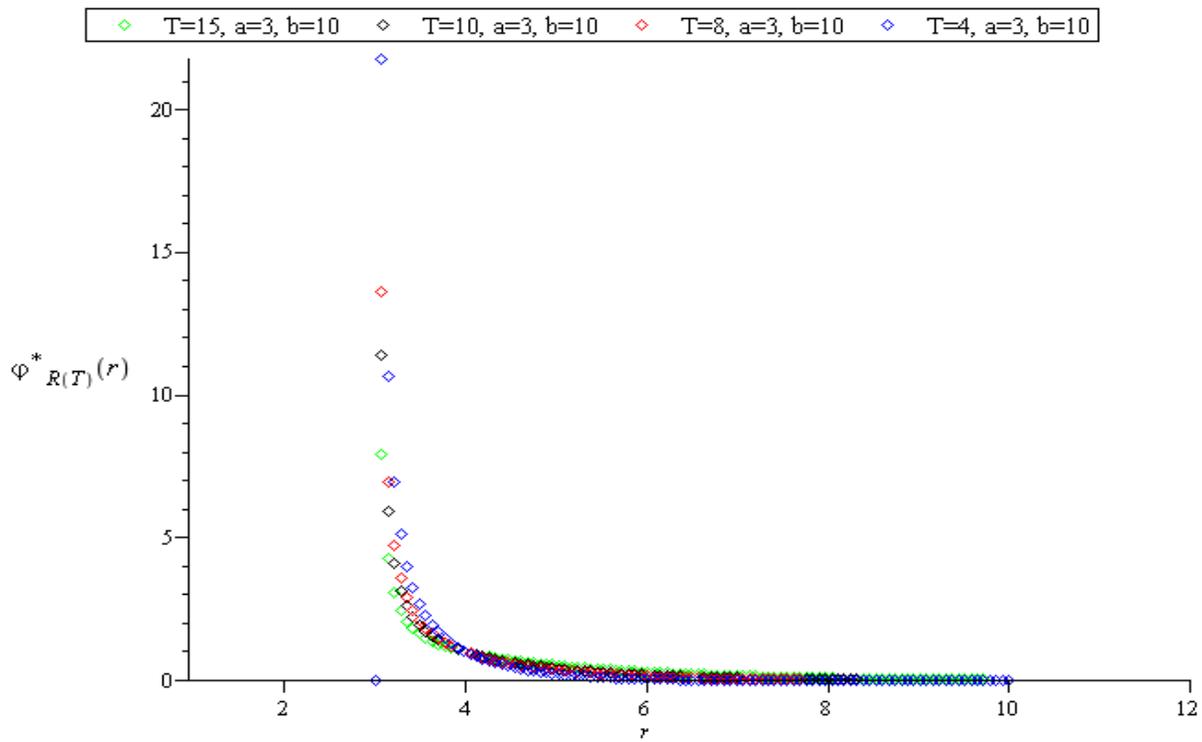

**Figure 11:** The reversed hazard rate function of TDDWPR.

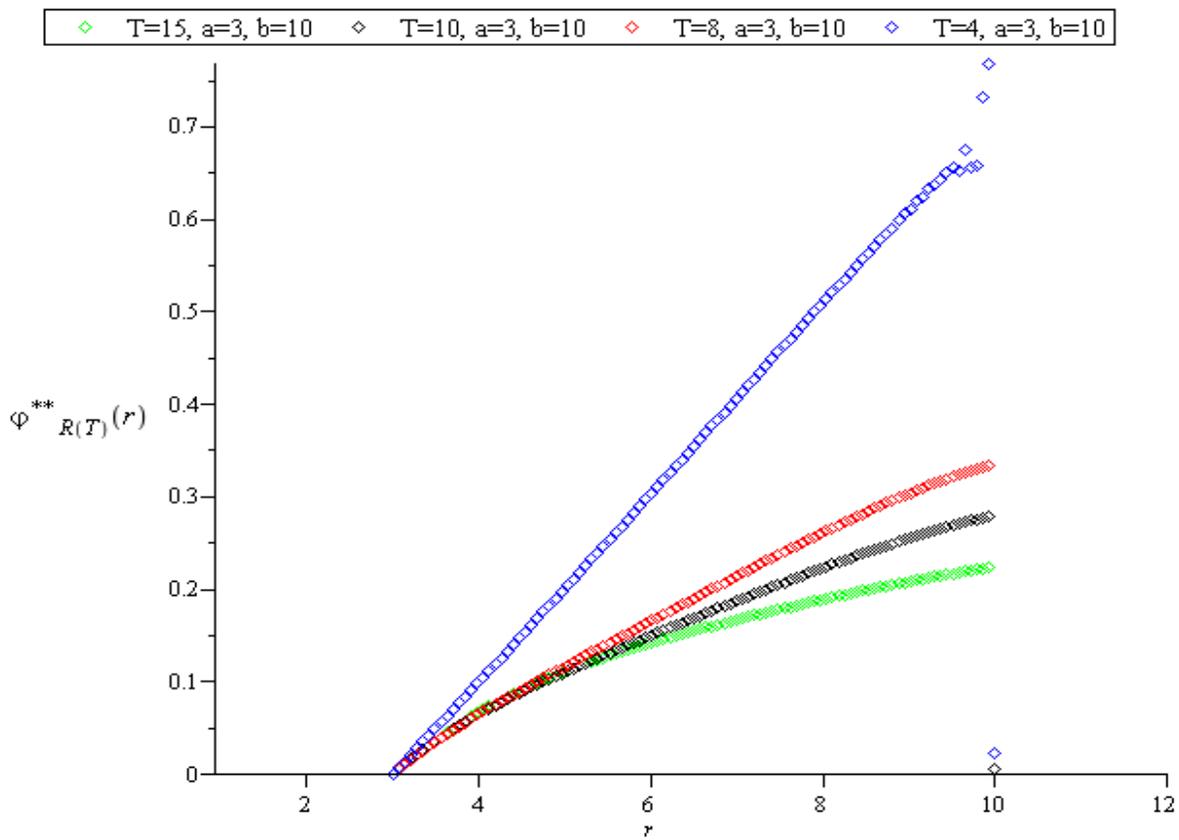

**Figure 12:** The second rate function of TDDWPR.

## 4.2 Moments of TDDWPR

By the same method in part (3.3) we can compute the value of the moments about the origin for a given values of $T$ as in the following Table 2. Let $\tilde{\mu}'_q$ donate to the $q^{th}$ moment for TDDWPR and compute from the following series:

$$\tilde{\mu}'_q = \sum_{q=0}^{\infty} \left[ \frac{r^q \sum_{k=1}^{\infty}\left[\left(\frac{8}{(2k-1)^2 \pi^2} + \frac{8T}{(r+1)^2}\right)\exp\left[-\frac{(2k-1)^2 \pi^2 T}{2(r+1)^2}\right] - \left(\frac{8}{(2k-1)^2 \pi^2} + \frac{8T}{r^2}\right)\exp\left[-\frac{(2k-1)^2 \pi^2 T}{2r^2}\right]\right]}{\sum_{k=1}^{\infty}\left[\left(\frac{8}{(2k-1)^2 \pi^2} + \frac{8T}{(b+1)^2}\right)\exp\left[-\frac{(2k-1)^2 \pi^2 T}{2(b+1)^2}\right]\right] - \sum_{k=1}^{\infty}\left[\left(\frac{8}{(2k-1)^2 \pi^2} + \frac{8T}{(a+1)^2}\right)\exp\left[-\frac{(2k-1)^2 \pi^2 T}{2(a+1)^2}\right]\right]} \right],$$

(23)

$r = a,\ldots,a+i,\ldots,b$, $i$ is a non negative real number.

**Table 2.** Moments, skewness and kurtosis for various values of $T$ and $a = 3, b = 10$ for TDDWPR

| $T$ | Raw moments | Central moments | Skewness | Kurtosis |
|---|---|---|---|---|
| 15 | $\tilde{\mu}'_1 = 3.791163482$<br>$\tilde{\mu}'_2 = 18.74262099$<br>$\tilde{\mu}'_3 = 95.95288124$<br>$\tilde{\mu}'_4 = 505.5189570$ | $\tilde{\mu}_2 = 4.36970044$<br>$\tilde{\mu}_3 = 1.766750584 \times 10^6$<br>$\tilde{\mu}_4 = 47.0013975$ | $1.934183818 \times 10^5$ | $-0.538456924$ |
| 20 | $\tilde{\mu}'_1 = 2.999409196$<br>$\tilde{\mu}'_2 = 15.67672675$<br>$\tilde{\mu}'_3 = 83.92201781$<br>$\tilde{\mu}'_4 = 458.2265663$ | $\tilde{\mu}_2 = 6.680271225$<br>$\tilde{\mu}_3 = 1.182052468 \times 10^6$<br>$\tilde{\mu}_4 = 54.7618919$ | $68461.37397$ | $-1.772870907$ |
| 25 | $\tilde{\mu}'_1 = 2.371112258$<br>$\tilde{\mu}'_2 = 12.82990901$<br>$\tilde{\mu}'_3 = 70.64012046$<br>$\tilde{\mu}'_4 = 394.5987442$ | $\tilde{\mu}_2 = 7.207735670$<br>$\tilde{\mu}_3 = 7.049715394 \times 10^5$<br>$\tilde{\mu}_4 = 62.58145720$ | $36431.18580$ | $-1.795385827$ |

## 4.3. Order Statistics of TDDWPR

If $\tilde{R}_{1:m} \leq \tilde{R}_{2:m} \leq ... \leq \tilde{R}_{m:m}$ be the order statistics of a random sample $\tilde{R}_1, \tilde{R}_2, ..., \tilde{R}_m$ from the TDDWPR, then the pmf of the $j^{th}$ order statistic $\tilde{R}_{j:m}$ is,

$$q_{\tilde{R}(T)(j:m)} = \frac{m!}{(j-1)!(m-j)!} \left(Q_{\tilde{R}(T)}(r)\right)^{j-1} \left(1 - Q_{\tilde{R}(T)}(r)\right)^{m-j} q_{\tilde{R}(T)}(r),$$

$$= \frac{m!}{(j-1)!(m-j)!}$$

$$\times \left( \frac{\sum_{k=1}^{\infty}\left[\left(\frac{8}{(2k-1)^2\pi^2} + \frac{8T}{(r+1)^2}\right)\exp\left[-\frac{(2k-1)^2\pi^2 T}{2(r+1)^2}\right] - \left(\frac{8}{(2k-1)^2\pi^2} + \frac{8T}{(a+1)^2}\right)\exp\left[-\frac{(2k-1)^2\pi^2 T}{2(a+1)^2}\right]\right]}{\sum_{k=1}^{\infty}\left[\left(\frac{8}{(2k-1)^2\pi^2} + \frac{8T}{(b+1)^2}\right)\exp\left[-\frac{(2k-1)^2\pi^2 T}{2(b+1)^2}\right] - \left(\frac{8}{(2k-1)^2\pi^2} + \frac{8T}{(a+1)^2}\right)\exp\left[-\frac{(2k-1)^2\pi^2 T}{2(a+1)^2}\right]\right]} \right)^{j-1}$$

$$\times \left(1 - \frac{\sum_{k=1}^{\infty}\left[\left(\frac{8}{(2k-1)^2\pi^2} + \frac{8T}{(r+1)^2}\right)\exp\left[-\frac{(2k-1)^2\pi^2 T}{2(r+1)^2}\right] - \left(\frac{8}{(2k-1)^2\pi^2} + \frac{8T}{(a+1)^2}\right)\exp\left[-\frac{(2k-1)^2\pi^2 T}{2(a+1)^2}\right]\right]}{\sum_{k=1}^{\infty}\left[\left(\frac{8}{(2k-1)^2\pi^2} + \frac{8T}{(b+1)^2}\right)\exp\left[-\frac{(2k-1)^2\pi^2 T}{2(b+1)^2}\right] - \left(\frac{8}{(2k-1)^2\pi^2} + \frac{8T}{(a+1)^2}\right)\exp\left[-\frac{(2k-1)^2\pi^2 T}{2(a+1)^2}\right]\right]}\right)^{m-j}$$

$$\times \left( \frac{\sum_{k=1}^{\infty}\left[\left(\frac{8}{(2k-1)^2\pi^2} + \frac{8T}{(r+1)^2}\right)\exp\left[-\frac{(2k-1)^2\pi^2 T}{2(r+1)^2}\right] - \left(\frac{8}{(2k-1)^2\pi^2} + \frac{8T}{r^2}\right)\exp\left[-\frac{(2k-1)^2\pi^2 T}{2r^2}\right]\right]}{\sum_{k=1}^{\infty}\left[\left(\frac{8}{(2k-1)^2\pi^2} + \frac{8T}{(b+1)^2}\right)\exp\left[-\frac{(2k-1)^2\pi^2 T}{2(b+1)^2}\right]\right] - \sum_{k=1}^{\infty}\left[\left(\frac{8}{(2k-1)^2\pi^2} + \frac{8T}{(a+1)^2}\right)\exp\left[-\frac{(2k-1)^2\pi^2 T}{2(a+1)^2}\right]\right]} \right).$$

(24)

The cdf of $\tilde{R}_{j:m}$ is also given by,

$$Q_{\tilde{R}(T)(j:m)}(r) = \sum_{i=j}^{m} \binom{m}{j} \left(Q_{\tilde{R}(T)}(r)\right)^{i} \left(1 - Q_{\tilde{R}(T)}(r)\right)^{m-i}$$

$$= \sum_{i=j}^{m} \binom{m}{i} \left( \frac{\sum_{k=1}^{\infty}\left[\left(\frac{8}{(2k-1)^2\pi^2} + \frac{8T}{(r+1)^2}\right)\exp\left[-\frac{(2k-1)^2\pi^2 T}{2(r+1)^2}\right] - \left(\frac{8}{(2k-1)^2\pi^2} + \frac{8T}{(a+1)^2}\right)\exp\left[-\frac{(2k-1)^2\pi^2 T}{2(a+1)^2}\right]\right]}{\sum_{k=1}^{\infty}\left[\left(\frac{8}{(2k-1)^2\pi^2} + \frac{8T}{(b+1)^2}\right)\exp\left[-\frac{(2k-1)^2\pi^2 T}{2(b+1)^2}\right] - \left(\frac{8}{(2k-1)^2\pi^2} + \frac{8T}{(a+1)^2}\right)\exp\left[-\frac{(2k-1)^2\pi^2 T}{2(a+1)^2}\right]\right]} \right)^{i}$$

$$\times \left(1 - \frac{\sum_{k=1}^{\infty}\left[\left(\frac{8}{(2k-1)^2\pi^2} + \frac{8T}{(r+1)^2}\right)\exp\left[-\frac{(2k-1)^2\pi^2 T}{2(r+1)^2}\right] - \left(\frac{8}{(2k-1)^2\pi^2} + \frac{8T}{(a+1)^2}\right)\exp\left[-\frac{(2k-1)^2\pi^2 T}{2(a+1)^2}\right]\right]}{\sum_{k=1}^{\infty}\left[\left(\frac{8}{(2k-1)^2\pi^2} + \frac{8T}{(b+1)^2}\right)\exp\left[-\frac{(2k-1)^2\pi^2 T}{2(b+1)^2}\right] - \left(\frac{8}{(2k-1)^2\pi^2} + \frac{8T}{(a+1)^2}\right)\exp\left[-\frac{(2k-1)^2\pi^2 T}{2(a+1)^2}\right]\right]}\right)^{m-i}.$$

(25)

The pmf of the minimum difference between the prices at knowing time $T$ is given by:

$$q_{\tilde{R}(T)(1:m)} = m \times \left(1 - \frac{\sum_{k=1}^{\infty}\left[\left(\frac{8}{(2k-1)^2\pi^2} + \frac{8T}{(r+1)^2}\right)\exp\left[-\frac{(2k-1)^2\pi^2 T}{2(r+1)^2}\right] - \left(\frac{8}{(2k-1)^2\pi^2} + \frac{8T}{(a+1)^2}\right)\exp\left[-\frac{(2k-1)^2\pi^2 T}{2(a+1)^2}\right]\right]}{\sum_{k=1}^{\infty}\left[\left(\frac{8}{(2k-1)^2\pi^2} + \frac{8T}{(b+1)^2}\right)\exp\left[-\frac{(2k-1)^2\pi^2 T}{2(b+1)^2}\right] - \left(\frac{8}{(2k-1)^2\pi^2} + \frac{8T}{(a+1)^2}\right)\exp\left[-\frac{(2k-1)^2\pi^2 T}{2(a+1)^2}\right]\right]}\right)^{m-1}$$

$$\times \left(\frac{\sum_{k=1}^{\infty}\left[\left(\frac{8}{(2k-1)^2\pi^2} + \frac{8T}{(r+1)^2}\right)\exp\left[-\frac{(2k-1)^2\pi^2 T}{2(r+1)^2}\right] - \left(\frac{8}{(2k-1)^2\pi^2} + \frac{8T}{r^2}\right)\exp\left[-\frac{(2k-1)^2\pi^2 T}{2r^2}\right]\right]}{\sum_{k=1}^{\infty}\left[\left(\frac{8}{(2k-1)^2\pi^2} + \frac{8T}{(b+1)^2}\right)\exp\left[-\frac{(2k-1)^2\pi^2 T}{2(b+1)^2}\right]\right] - \sum_{k=1}^{\infty}\left[\left(\frac{8}{(2k-1)^2\pi^2} + \frac{8T}{(a+1)^2}\right)\exp\left[-\frac{(2k-1)^2\pi^2 T}{2(a+1)^2}\right]\right]}\right).$$

(26)

and the pmf of the maximum difference is,

$$q_{\tilde{R}(T)(m:m)} = m \left(\frac{\sum_{k=1}^{\infty}\left[\left(\frac{8}{(2k-1)^2\pi^2} + \frac{8T}{(r+1)^2}\right)\exp\left[-\frac{(2k-1)^2\pi^2 T}{2(r+1)^2}\right] - \left(\frac{8}{(2k-1)^2\pi^2} + \frac{8T}{(a+1)^2}\right)\exp\left[-\frac{(2k-1)^2\pi^2 T}{2(a+1)^2}\right]\right]}{\sum_{k=1}^{\infty}\left[\left(\frac{8}{(2k-1)^2\pi^2} + \frac{8T}{(b+1)^2}\right)\exp\left[-\frac{(2k-1)^2\pi^2 T}{2(b+1)^2}\right] - \left(\frac{8}{(2k-1)^2\pi^2} + \frac{8T}{(a+1)^2}\right)\exp\left[-\frac{(2k-1)^2\pi^2 T}{2(a+1)^2}\right]\right]}\right)^{m-1}$$

$$\times \left(\frac{\sum_{k=1}^{\infty}\left[\left(\frac{8}{(2k-1)^2\pi^2} + \frac{8T}{(r+1)^2}\right)\exp\left[-\frac{(2k-1)^2\pi^2 T}{2(r+1)^2}\right] - \left(\frac{8}{(2k-1)^2\pi^2} + \frac{8T}{r^2}\right)\exp\left[-\frac{(2k-1)^2\pi^2 T}{2r^2}\right]\right]}{\sum_{k=1}^{\infty}\left[\left(\frac{8}{(2k-1)^2\pi^2} + \frac{8T}{(b+1)^2}\right)\exp\left[-\frac{(2k-1)^2\pi^2 T}{2(b+1)^2}\right]\right] - \sum_{k=1}^{\infty}\left[\left(\frac{8}{(2k-1)^2\pi^2} + \frac{8T}{(a+1)^2}\right)\exp\left[-\frac{(2k-1)^2\pi^2 T}{2(a+1)^2}\right]\right]}\right).$$

(27)

### 4.4. Stress-strength parameter of TDDWPR

If the double truncated random variable $\tilde{R}(T)$ (which describes the changing of stock price and take the values $r = a,,...,a+i,...,b$, ($i$ is a non negative real number)) is consider as the strength of a component which is subjected to a random stress $\tilde{Z}$, then the stress-strength model is defined as,

$$P(\tilde{R}(T) > \tilde{Z}) = \sum_{r=0}^{\infty} q_{\tilde{R}(T)}(r) Q_{\tilde{R}(T)}(r)$$

By using (19), (20) and if $\tilde{R}(T)$ has $TDDWPR(T_3)$ and $\tilde{Z}$ has $TDDWPR(T_4)$ then we get the stress-strength parameter by,

$$P(\tilde{R}(T) > \tilde{Z}) = \frac{1}{\xi_k^2} \sum_{r=0}^{\infty} \left( \sum_{k=1}^{\infty} \left[ \left( \frac{8}{(2k-1)^2 \pi^2} + \frac{8T}{(r+1)^2} \right) \exp\left[ -\frac{(2k-1)^2 \pi^2 T}{2(r+1)^2} \right] - \left( \frac{8}{(2k-1)^2 \pi^2} + \frac{8T}{r^2} \right) \exp\left[ -\frac{(2k-1)^2 \pi^2 T}{2r^2} \right] \right] \right.$$
$$\left. \times \sum_{k=1}^{\infty} \left[ \left( \frac{8}{(2k-1)^2 \pi^2} + \frac{8T}{(r+1)^2} \right) \exp\left[ -\frac{(2k-1)^2 \pi^2 T}{2(r+1)^2} \right] - \left( \frac{8}{(2k-1)^2 \pi^2} + \frac{8T}{(a+1)^2} \right) \exp\left[ -\frac{(2k-1)^2 \pi^2 T}{2(a+1)^2} \right] \right] \right),$$

where,

$$\xi_k = \sum_{k=1}^{\infty} \left[ \left( \frac{8}{(2k-1)^2 \pi^2} + \frac{8T}{(b+1)^2} \right) \exp\left[ -\frac{(2k-1)^2 \pi^2 T}{2(b+1)^2} \right] \right] - \sum_{k=1}^{\infty} \left[ \left( \frac{8}{(2k-1)^2 \pi^2} + \frac{8T}{(a+1)^2} \right) \exp\left[ -\frac{(2k-1)^2 \pi^2 T}{2(a+1)^2} \right] \right].$$

## 5. Application

Economists evaluate the stock price by specifying the points where the number of buyers are more than the sellers and their bids to buy for the highest possible price and stop the price at these points because of the desire of buyers to buy higher than this price or vice versa. Since, the differences between these points at given time intervals are random, we use DDWPR distribution. It is clear that, the data of this distribution is depend on the stochastic Weiner model, so these data are dependent on each other. Therefore, we consider here the same values of the data which taken in Withers and Nadarajah [23].

Due to the importance of the cumulative distributions of $R$ and $\tilde{R}$ over a certain period of time $T$, we obtained the values of the cdf and pmf of DDWPR and TDDWPR for $r, r = 1,2,...,5$ and the corresponding time periods $T, T = 25, 50, 75, 100, 200$ are given in Table 3. This is because the non-accumulative shareholders lose their share of profits for any period during which the members of the board of directors do not announce a dividend.

Where the earnings per share over $T$ is the share of the stock or the share price difference of the income period $T$, so there is a relationship between that difference and the length of $T$. The advantage of DDWPR and its truncated version is appear from the numerical calculations. We notice that their exist a positive relationship between the values of the difference and $T$ where the values of the CDF of DDWPR and TDDWPR are converge to 1 when $R \geq T$ and $\tilde{R} \geq T$, respectively, see figures 2 and 8. Also, we notice that when the length of the time is increasing the efficiency of this distribution is increasing. We conclude that this distribution is a one of the important lifetime distribution which has a wide applications in the real life.

**Table 3:** The pmf and the cdf of DDWPR and TDDWPR.

| $T$ | $r$ | $f_{R(T)}(r)$ | $F_{R(T)}(r)$ | $a=3,\ b=10$ | |
|---|---|---|---|---|---|
| | | | | $q_{\tilde{R}(T)}(r)$ | $Q_{\tilde{R}(T)}(r)$ |
| 25 | 1 | $2.112541439 \times 10^{-12}$ | $2.112541439 \times 10^{-12}$ | $2.3895748 \times 10^{-12}$ | 0 |
| | 2 | 0.00002601247475 | 0.00002601247686 | 0.00002942368534 | 0 |
| | 3 | 0.005984818378 | 0.006010830857 | 0.006769652425 | 0 |
| | 4 | 0.05767745832 | 0.06368828918 | 0.06524113530 | 0.065241135 |
| | 5 | 0.1438681320 | 0.2075564212 | 0.1627346374 | 0.227975772 |
| | 9 | 0.1033830971 | 0.8197700160 | 0.1169404967 | 0.9204735202 |
| | 10 | 0.0703066431 | 0.8900766592 | 0.07952647965 | 1.000000000 |
| 50 | 1 | $1.742600764 \times 10^{-25}$ | $1.742600764 \times 10^{-25}$ | $3.245690482 \times 10^{-25}$ | 0 |
| | 2 | $5.771834028 \times 10^{-11}$ | $5.771834028 \times 10^{-11}$ | $1.075036070 \times 10^{-10}$ | 0 |
| | 3 | 0.00005262943830 | 0.00005263001549 | 0.000009802524507 | 0 |
| | 4 | 0.0008730186293 | 0.0008782816312 | 0.001626045572 | 0.0016260455 |
| | 5 | 0.01179173697 | 0.01267001860 | 0.02196276350 | 0.0235888090 |
| | 9 | 0.1348648062 | 0.4090517634 | 0.2511931745 | 0.7618717726 |
| | 10 | 0.1278502787 | 0.5369020420 | 0.2381282284 | 1.00000000 |

**Table 3:** continued

| T | r | $f_{R(T)}(r)$ | $F_{R(T)}(r)$ | $a=3,\ b=10$ | |
|---|---|---|---|---|---|
| | | | | $q_{\tilde{R}(T)}(r)$ | $Q_{\tilde{R}(T)}(r)$ |
| | 1 | $1.08384351 \times 10^{-38}$ | $1.083843510 \times 10^{-38}$ | $3.988800472 \times 10^{-38}$ | 0 |
| | 2 | $9.718996748 \times 10^{-17}$ | $9.718996752 \times 10^{-17}$ | $3.576820681 \times 10^{-16}$ | 0 |
| | 3 | $3.527451308 \times 10^{-9}$ | $3.527451406 \times 10^{-9}$ | $1.298185514 \times 10^{-8}$ | 0 |
| 75 | 4 | 0.000009365538985 | 0.000009369066440 | 0.00003446739864 | 0.0022283382 |
| | 5 | 0.0005961222298 | 0.0006054912962 | 0.002193870803 | 0.0000344674 |
| | 9 | 0.08326008312 | 0.1688484851 | 0.3064167989 | 0.6214023490 |
| | 10 | 0.1028731848 | 0.2717216698 | 0.3785976521 | 1.00000000 |
| | 1 | $6.000188914 \times 10^{-52}$ | $6.000188914 \times 10^{-52}$ | $4.753316306 \times 10^{-51}$ | 0 |
| | 2 | $1.459057072 \times 10^{-22}$ | $1.459057071 \times 10^{-22}$ | $1.155856902 \times 10^{-21}$ | 0 |
| | 3 | $2.112541439 \times 10^{-12}$ | $2.112541439 \times 10^{-12}$ | $1.673543586 \times 10^{-11}$ | 0 |
| 100 | 4 | $8.955187063 \times 10^{-8}$ | $8.955398320 \times 10^{-8}$ | $7.094249415 \times 10^{-7}$ | $7.094249 \times 10^{-7}$ |
| | 5 | 0.00002592292288 | 0.00002601247686 | 0.0002053599542 | 0.0002060693 |
| | 9 | 0.03938290275 | 0.06368828918 | 0.3119891665 | 0.5045350868 |
| | 10 | 0.06254334652 | 0.1262316357 | 0.4954649147 | 1.00000000 |

## 6. Concluding remarks

In this paper, we introduce the discrete distribution for a Wiener process range (DDWPR). This distribution is the best for the stock price to obtain cumulative shareholders and do not lose their share of profits for any period of time. We provided a mathematical treatment to study the basic statistical properties of this discrete model including reliability properties, moments, stress-strength parameter, quantile's and its random number generation and order statistics. Furthermore, we study the truncated version of this distribution (TDDWPR) and its prosperities as a one important distribution for the real life problem. A data set is analyzed to clarify the effectiveness of DDWPR and TDDWPR. Some issues for future research may be considering to study the corresponding step-stress model for our distribution.